\documentclass[leqno,12pt]{amsart} 

\usepackage{amssymb, amsmath}

\theoremstyle{plain} 
\newtheorem{thm}{Theorem}[section]
\newtheorem{lem}[thm]{Lemma}
\newtheorem{prop}[thm]{Proposition}

\theoremstyle{definition}

\newtheorem{ex}[thm]{Example}
\newtheorem*{ackn}{Acknowledgement}

\newcommand{\skal}[2]{\langle #1,#2\rangle}

\begin{document}

\title[On leafwise conformal diffeomorphisms]{On leafwise conformal diffeomorphisms}
\author{Kamil Niedzia\l omski}
\address{
Department of Mathematics and Computer Science \endgraf
University of \L\'{o}d\'{z} \endgraf
ul. Banacha 22, 90-238 \L\'{o}d\'{z} \endgraf
Poland
}
\email{kamiln@math.uni.lodz.pl}

\maketitle

\begin{abstract} For every diffeomorphism $\varphi:M\to N$ between $3$--dimensional Riemannian manifolds $M$ and $N$ there are in general locally two $2$--dimensional distributions $D_{\pm}$ such that $\varphi$ is conformal on both of them. We state necessary and sufficient conditions for a distribution to be one of $D_{\pm}$. These are algebraic conditions expressed in terms of the self-adjoint and positive definite operator $(\varphi_{\ast})^*\varphi_{\ast}$. We investigate integrability condition of $D_+$ and $D_-$. We also show that it is possible to choose coordinate systems in which leafwise conformal diffeomorphism is holomorphic on leaves.
\end{abstract}

\section{Introduction} Let $\varphi:M\to N$ be a diffeomorphism between $3$--dimensional Riemannian manifolds $(M,g)$ and $(N,h)$. Fix $x\in M$ and let $(\varphi_{\ast x})^*:T_{\varphi(x)}N\to T_xM$ denotes the operator adjoint to $\varphi_{\ast x}:T_xM\to T_{\varphi(x)}N$. Then $S_x=(\varphi_{\ast x})^*\varphi_{\ast x}$ is a self--adjoint and positive definite operator. Let $0<\lambda_1(x)\leq\lambda_2(x)\leq\lambda_3(x)$ be the eigenvalues of $S_x$. 

Preimage $E(x)=\varphi_{\ast x}^{-1}(\mathbb{S}^2)$ of the unit sphere is an ellipsoid with principial semi--axes $1/\sqrt{\lambda_i(x)}$, $i=1,2,3$. Therefore, if the eigenvalues $\lambda_i(x), i=1,2,3,$ are distinct there are two $2$--dimensional subspaces $D_{+}(x)$ and $D_-(x)$ of $T_xM$ intersecting $E(x)$ along spheres. Thus locally we get two smooth distributions $D_+$ and $D_-$. By the definition of $D_{\pm}$ we see that $\varphi$ is conformal on each of them (see Lemma \ref{l1}).

In this article we describe $D_+$ and $D_-$ and study the problem of integrability of these distributions. We show that integrability of one of the distributions $D_{\pm}$ does not imply integrability of the other one.

Conformality of diffeomorphisms on distributions of codimension one was studied by S. Tanno in \cite{t1} and \cite{t2}. However, majority of results in \cite{t1} and \cite{t2} is obtained under the assumption that a given diffeomorphism $\varphi$ maps vectors normal to a distribution $D$ to vectors normal to the image $\varphi_{\ast}(D)$. Therefore $\varphi$ cannot have distinct eigenvalues. Moreover, in \cite{kn} the author showed that under some assumptions on a diffeomorphisms $\varphi$ and the dimension of $M$  there are no distributions of `small' codimension on which $\varphi$ is conformal. In particular, assuming $\dim M>3$ there are no codimension one foliations such that a diffoemorphism $\varphi:M\to N$, for which $S$ has distinct eigenvalues, is conformal on the leaves.

The paper is organized as follows. In section 2 we obtain preliminary results concerning some operators defined for $1$--forms. Next, we state necessary and sufficient conditions for a diffeomorphism between $3$--dimensional Riemannian manifolds to be conformal on a given distribution, that is we obtain conditions for a distribution to be one of $D_{\pm}$ (Theorem \ref{t1}). Examples are given. In the following sections we focus on the integrability condition of $D_+$ and $D_-$ (Theorem \ref{t2}, Propositions \ref{p2} and \ref{p3}). The last part of this article is devoted to local description of leafwise conformal diffeomorphism. We show that it is possible to choose appropriate coordinate systems in which given leafwise conformal diffeomorphism is holomorphic on leaves (Theorem \ref{t3}).

\section{Notations and preliminary results} Let $(M,g)$, $(N,h)$ be $3$--dimensional oriented and connected Riemannian manifolds and let $\varphi:(M,g)\to (N,h)$ be a diffeomorphism. We say that $\varphi$ is {\it leafwise conformal} if there exists a $2$--dimensional foliation $\mathcal{F}$ on $M$ such that $\varphi:L\to \varphi(L)$ is conformal for every leaf $L\in\mathcal{F}$. In that case we also say that $\varphi$ is $\mathcal{F}$--{\it conformal}. $\varphi$ is {\it locally leafwise conformal} if every point $x\in M$ has a neighbourhood $U$ such that $\varphi:U\to \varphi(U)$ is leafwise conformal.

Let $\lambda_1,\lambda_2,\lambda_3$ be the eigenvalues of the operator $S=(\varphi_{\ast})^*\varphi_{\ast}:TM\to TM$ and $\xi_1,\xi_2,\xi_3$ be the corresponding unit eigenvectors. Assume $\lambda_1<\lambda_2<\lambda_3$. Let $\eta_1,\eta_2,\eta_3$ be the basis dual to $\xi_1,\xi_2,\xi_3$. Locally we may choose above bases to be smooth. Define
\begin{equation} \label{eq1}
\omega_{\pm}=\frac{\sqrt{\lambda_2-\lambda_1}}{\sqrt{\lambda_3-\lambda_1}}\eta_1\pm\frac{\sqrt{\lambda_3-\lambda_2}}{\sqrt{\lambda_3-\lambda_1}}\eta_3.
\end{equation}
Condsider the distributions $D_{\pm}=\ker\omega_{\pm}$. We have
\begin{lem} \label{l1} A diffeomorphism $\varphi$ is (locally) conformal on a $2$--dimensional distribution $D$ if and only if $D=D_+$ or $D=D_-$ (locally). Moreover the coefficient of conformality is $\lambda_2$.
\end{lem}
\begin{proof} It is easy to check that $\varphi$ is conformal on $D_+$ and $D_-$ with coefficient of conformality $\lambda_2$. Suppose there exists a distribution $D$ such that $\varphi$ is conformal on $D$. Fix $x\in M$ and consider the set $E(x)=d\varphi^{-1}(x)(\mathbb{S}^2)$, where $\mathbb{S}^2\subset T_{\varphi(x)}N$ is the unit sphere. Then $E(x)$ is an ellipsoid with principial semi--axes $1/\sqrt{\lambda_i(x)}$, $i=1,2,3$. The subspaces $D_+(x)$ and $D_-(x)$ intersect $E(x)$ along circles and these are the only subspaces with this property, see \cite{kb} or \cite{kn}. Thus by conformality of $\varphi$ on $D$ we get that $D(x)=D_+(x)$ or $D(x)=D_-(x)$. Since $M$ is connected, $D$ is smooth and $D_+(x)\not=D_-(x)$ for all $x\in M$, we obtain $D=D_+$ or $D=D_-$ (locally).
\end{proof} 
Let $x\in M$, $p=0,1,2,3$ and $\ast:\Lambda^pT^*_xM\to\Lambda^{3-p}T^*_xM$ be the Hodge operator. Let $\iota(\omega)\eta=\omega\wedge\eta$ for $\omega,\eta\in \Lambda^pT^*_xM$. For $\omega,\eta\in T^*_xM$ define $(\omega\odot\eta)_x:T^*_xM\to T^*_xM$ by
\[
(\omega\odot\eta)_x\alpha=\skal{\omega}{\alpha}\eta+\skal{\eta}{\alpha}\omega,\quad\alpha\in T^*_xM,
\]
where $\skal{\cdot}{\cdot}$ is the inner product in $T^*_xM$ induced from Riemannian metric $g$.
Moreover for $\theta\in[0,2\pi)$ and $\omega\in T^*_xM$, $|\omega|=1$, put
\[
{\rm Rot}_x(\theta,\omega)={\rm Id}_{T^*_xM}+\sin\theta(\ast\iota(\omega))+(1-\cos\theta)(\ast\iota(\omega))^2:T^*_xM\to T^*_xM.
\]
Then ${\rm Rot}_x(\theta,\omega)$ is an operator of rotation around $\omega$ of an angle $\theta$, for details see \cite{ra}. For simplicty we will write ${\rm Rot}_x(\omega)$ instead of ${\rm Rot}_x(\pi/2,\omega)$. 
\begin{lem} \label{l2} Let $0\leq \theta,\theta_1,\theta_2<2\pi$, $\omega,\eta\in T^*_xM$ and $|\omega|=1$. The operator ${\rm Rot}_x(\theta,\omega)$ has the following properties
\begin{enumerate}
\item[$(1)$] ${\rm Rot}_x(\theta_1,\omega)\circ {\rm Rot}_x(\theta_2,\omega)={\rm Rot}_x(\theta_1+\theta_2 \mod 2\pi,\omega)$.
\item[$(2)$] If $\skal{\omega}{\eta}=0$ then $\skal{\omega}{{\rm Rot}_x(\theta,\omega)\eta}=0$.
\item[$(3)$] If $\skal{\omega}{\eta}=0$ then $\skal{{\rm Rot}_x(\omega)\eta}{\eta}=0$ and $\eta-{\rm Rot}_x(\omega)\eta=\sqrt{2}{\rm Rot}_x(-\frac{\pi}{4},\omega)\eta$.
\end{enumerate}
\end{lem}
\begin{proof} Easy computations left to the reader.
\end{proof}

The operator $S_x:T_xM\to T_xM$ can be considered as an operator $S_x:T^*_xM\to T^*_xM$ by the rule $(S_x\eta)X=\eta(SX)$, $X\in T_xM$. Then $S$ is a self--adjoint and positive definite operator with eigenvalues $\lambda_i$ and corresponding eigenvectors $\eta_i$, $i=1,2,3$. Let $[T_1,T_2]=T_1T_2-T_2T_1:T^*_xM\to T^*_xM$ be the comutator of operators $T_1,T_2:T^*_xM\to T^*_xM$. We define
\begin{align}
B_x(\omega) &=[S_x,\ast\iota(\omega)]:T^*_xM\to T^*_xM, \label{defb} \\
A_x(\omega) &=[S_x,{\rm Rot}_x(\omega)]:T^*_xM\to T^*_xM. \label{defa}
\end{align}

We have a technical result
\begin{lem} \label{l3} Let $\omega\in T^*_xM$. Then there exist $\eta,\sigma\in T^*_xM$ such that $\omega,\eta,\sigma$ are orthogonal and
\begin{equation} \label{l3e}
S_x\eta=\frac{1}{|\eta|^2}\eta+\skal{S_x\omega}{\eta}\omega,\quad 
S_x\sigma=\frac{1}{|\sigma|^2}\sigma+\skal{S_x\omega}{\sigma}\omega.
\end{equation}
\end{lem}
\begin{proof} Let $\omega=\sum_ia_i\eta_i$. If $\omega=\eta_i$ for some $i=1,2,3$, then it sufficies to put $\eta=(1/\sqrt[3]{\lambda_j})\eta_j$ and $\sigma=(1/\sqrt[3]{\lambda_k})\eta_k$, where $(i,j,k)$ is a permutation of the set $\{1,2,3\}$. 

Suppose now $\omega\not=\eta_i$ for all $i=1,2,3$. Let $C>0$ be such that $\sum_ia_i^2/(\lambda_i-C)=0$ and put $\eta=\sum_i(a_i/(\lambda_i-C))\eta_i$. Then $\skal{\omega}{\eta}=0$ and $S_x\eta=C\eta+\omega$. It sufficies to multiply $\eta$ by $1/\sqrt{C}|\eta|$. Let $\sigma={\rm Rot}_x(\omega)\eta$. By Lemma \ref{l2} $\omega,\eta,\sigma$ are orthogonal. Moreover, $\skal{S_x\sigma}{\eta}=0$ and $\skal{S_x\sigma}{\sigma}>0$, thus multiplying $\sigma$ by an appropriate factor we get $S_x\sigma=\frac{1}{|\sigma|^2}\sigma+\skal{S\omega}{\sigma}\omega$.
\end{proof}

\section{Conformality on distribution} 
Let $(M,g)$, $(N,h)$ be $3$--dimensional oriented and connected Riemannian manifolds and let $\varphi:(M,g)\to (N,h)$ be a diffeomorphism. Consider notations from the previous section.
\begin{thm} \label{t1} Let $D=\ker\omega$ be a $2$--dimensional distribution on an open subset $U\subset M$, where $\omega$ is a unit $1$--form on $U$. Assume the operator $S$ has distinct eigenvalues $\lambda_1<\lambda_2<\lambda_3$ and the corresponding unit eigenvectors $\eta_1, \eta_2, \eta_2$ are smooth on $U$. Then the following conditions are equivalent
\begin{enumerate}
\item[$(1)$] $\varphi$ is conformal on $D$, 
\item[$(2)$] $B(\omega)=\mu(\omega\odot\eta_2)$ for some smooth and nowhere vanishing function $\mu$ on $U$, 
\item[$(3)$] $A(\omega)^3=0$.
\end{enumerate}
Moreover, if $(2)$ holds then $\mu$ is equal to
\begin{equation} \label{mu}
\mu=\sqrt{\lambda_2-\lambda_1}\sqrt{\lambda_3-\lambda_2}\quad\textrm{or}\quad \mu=-\sqrt{\lambda_2-\lambda_1}\sqrt{\lambda_3-\lambda_2}.
\end{equation}
\end{thm}
\begin{proof}
$(1)\Rightarrow (2)$ The $1$--form $\omega$ is given by (\ref{eq1}) with sign $+$ or $-$ in place of $\pm$. Therefore with respect to the basis $\{\eta_1,\eta_2,\eta_3\}$ the operator $\ast\iota(\omega)$ is  represented by the matrix
\[
\ast\iota(\omega)=\left[\begin{array}{ccc}
0 & \pm\frac{\sqrt{\lambda_3-\lambda_2}}{\sqrt{\lambda_3-\lambda_1}} & 0 \\
\mp\frac{\sqrt{\lambda_3-\lambda_2}}{\sqrt{\lambda_3-\lambda_1}} & 0 & -\frac{\sqrt{\lambda_2-\lambda_1}}{\sqrt{\lambda_3-\lambda_1}} \\
0 & \frac{\sqrt{\lambda_2-\lambda_1}}{\sqrt{\lambda_3-\lambda_1}} & 0
\end{array}\right].
\] 
Since $S$ is represented by a diagonal matrix ${\rm diag}(\lambda_1,\lambda_2,\lambda_3)$ easy computations lead to equality $B(\omega)=\mu (\omega\odot\eta_2)$, where $\mu=\pm\sqrt{\lambda_3-\lambda_2}\sqrt{\lambda_2-\lambda_1}$.

$(2)\Rightarrow (3)$ Since for any two $1$--forms $\alpha,\beta$ we have ${\rm Tr}(\alpha\odot\beta)=2\skal{\alpha}{\beta}$ then
\[
0={\rm Tr}B(\omega)=\mu{\rm Tr}(\omega\odot\eta_2)=2\mu\skal{\omega}{\eta_2}.
\]
Thus $\omega$ and $\eta_2$ are orthonormal. Let $\sigma$ be a $1$--form such that $\{\omega,\eta_2,\sigma\}$ is an oriented orthonormal basis. Then with respect to this basis $B(\omega)$ and $\ast\iota(\omega)$ are represented by matrices
\[
B(\omega)=\mu\left[\begin{array}{ccc} 0 & 1 & 0 \\ 1 & 0 & 0 \\ 0 & 0 & 0 \end{array}\right],\quad 
\ast\iota(\omega)=\left[\begin{array}{ccc} 0 & 0 & 0 \\ 0 & 0 & -1 \\ 0 & 1 & 0 \end{array}\right].
\]
Since $A(\omega)=B(\omega)+(\ast\iota(\omega))B(\omega)+B(\omega)(\ast\iota(\omega))$, we have 
\[
A(\omega)=\mu\left[\begin{array}{ccc} 0 & 1 & -1 \\ 1 & 0 & 0 \\ 1 & 0 & 0 \end{array}\right],
\]
thus $A(\omega)^3=0$. 

$(3)\Rightarrow (1)$ Suppose $\varphi$ is not conformal on $D=\ker\omega$. Then $\varphi$ is not conformal on $D(x)$ at some point $x\in M$. Consider a set $L=\varphi_{\ast}^{-1}(\mathbb{S}^2)\cap D(x)$, where $\mathbb{S}^2\subset T_{\varphi(x)}N$ is the unit sphere. Then $L$ is an ellipse but not a circle. Let $\eta$ and $\sigma$ be as in Lemma \ref{l3}. Then ${\rm Rot}(\omega)\eta=a_1\sigma$ and ${\rm Rot}(\omega)\sigma=a_2\eta$ where $a_1a_2=-1$. Put $C=1/|\eta|^2-1/|\sigma|^2$. Then by \eqref{l3e}
\[
A(\omega)\eta=a_1C\sigma+e_1\omega,\quad A(\omega)\sigma=a_2C\eta+e_2\omega,
\]
where 
\[
e_1=\skal{S\omega}{a_1\eta-\sigma},\quad e_2=\skal{S\omega}{-a_2\eta-\sigma}.
\]
We have $S\omega={\rm Const}\cdot\omega+\omega_0$, where $\skal{\omega}{\omega_0}=0$. Put $\tau={\rm Rot}(-\pi/4,\omega)\omega_0$. Then by Lemma \ref{l2} $\skal{\omega}{\tau}=0$. Hence, $\tau=b_{\eta}\eta+b_{\sigma}\sigma$ for some $b_{\eta},b_{\sigma}$. Moreover, ${\rm Rot}(\omega)\tau=b_{\eta}a_1\sigma-b_{\sigma}a_2\eta$. Hence $b_{\sigma}e_1+b_{\eta}e_2=0$. Therefore
\[
A(\omega)\tau=C(b_{\eta}a_1\sigma+b_{\sigma}a_2\eta).
\]
By Lemma \ref{l2} $A(\omega)\omega=\sqrt{2}\tau$. Hence by assumption $A(\omega)^3=0$ we have
\[
0=A(\omega)^3\omega=A(\omega)^2\tau=-C^2\tau+C(b_{\eta}a_1e_2+b_{\sigma}a_2e_1)\omega.
\]
Since $S$ has distinct eigenvalues, then $\omega_0\not=0$ and $\tau\not=0$. Thus by linear indepedance of $\tau$ and $\omega$ we have $C=0$. Therefore $|\eta|=|\sigma|$. Since $\skal{S\eta}{\eta}=\skal{S\sigma}{\sigma}=1$ it follows that $L$ is a circle. Contradiction.
\end{proof}
\begin{ex} \label{e1} Let $U=\{x=(x_1,x_2,x_3)\in\mathbb{R}^3: \cos(x_2+x_3)\not=0\}$. Define a map $\varphi:U\to\mathbb{R}^3$ between Euclidean spaces in the following way
\[
\varphi(x)=(-\cos x_2+\sqrt{2}\sin x_3,\sin x_2-\sqrt{2}\cos x_3,\sqrt{2}x_1+x_2),\quad x=(x_1,x_2,x_3).
\]
Then
\[
S(x)=\left[\begin{array}{ccc} 2 & \sqrt{2} & 0 \\ \sqrt{2} & 2 & \sqrt{2}\sin(x_2+x_3) \\ 0 & \sqrt{2}\sin(x_2+x_3) & 2 
\end{array}\right]
\]
and $\det S(x)=4\cos^2(x_2+x_3)$. Therefore $\varphi$ is a diffeomorphism on $U$. Moreover the eigenvalues of $S$ are 
\[
2-\sqrt{2+2\sin^2(x_2+x_3)},\quad 2,\quad 2+\sqrt{2+2\sin^2(x_2+x_3)}.
\]
Thus $S(x)$ has distinct eigenvalues for every $x\in U$. Put
\[
\omega_+=dx_2, \quad \omega_-=\frac{\sqrt{2}}{\sqrt{2+2\sin^2(x_2+x_3)}}dx_1+\frac{\sqrt{2}\sin(x_2+x_3)}{\sqrt{2+2\sin^2(x_2+x_3)}}dx_3.
\]
Then $A(\omega_+)$ equals to
\[
\left[\begin{array}{ccc}
0 & -\sqrt{2}(\sin(x_2+x_3)-1) & 0 \\ -\sqrt{2}(\sin(x_2+x_3)+1) & 0 & -\sqrt{2}(\sin(x_2+x_3)-1) \\
0 & \sqrt{2}(\sin(x_2+x_3)+1) & 0
\end{array}\right]
\]
and $A(\omega_-)$ to
\[
\left[\begin{array}{ccc}
4\frac{\sin(x_2+x_3)}{\sqrt{2+2\sin^2(x_2+x_3)}} & -\sqrt{2} & -2\frac{\cos^2(x_2+x_3)}{\sqrt{2+2\sin^2(x_2+x_3)}} \\
\sqrt{2} & 0 & \sqrt{2}\sin(x_2+x_3) \\
-2\frac{\cos^2(x_2+x_3)}{\sqrt{2+2\sin^2(x_2+x_3)}} & -\sqrt{2}\sin(x_2+x_3) & -4\frac{\sin(x_2+x_3)}{\sqrt{2+2\sin^2(x_2+x_3)}}
\end{array}\right].
\]
Therefore $A(\omega_+)^3=A(\omega_-)^3=0$. Thus $\varphi$ is conformal on distributions $D_+=\ker\omega_+$ and $D_-=\ker\omega_-$.
\end{ex}

\section{Integrability condition in terms of an orthonormal moving frame} Let $(M,g)$ be a $3$--dimensional oriented Riemannian manifold. Let $e_1,e_2,e_3$ be a local orthonormal basis on the open subset $U$ of $M$ and $\omega_1,\omega_2,\omega_3$ the dual basis of $1$--forms. Let $(N,h)$ be another $3$--dimensional oriented Riemannian manifold and $\varphi:M\to N$ be a diffeomorphism. Let $D=\ker\omega$ be a two dimensional distribution on $U$, where $\omega$ is a unit $1$--form on $M$. Suppose $S=(\varphi_{\ast})^*\varphi_{\ast}$ has distinct eigenvalues. Let $\mu$ be a smooth nowhere vanishing function on $U$. Let $\lambda$ be the middle eigenvalue of $S$ and $\eta_2$ the unit eivenvector corresponding to $\lambda$. Let
\[
\omega=\sum_i\alpha_i\omega_i,\quad \eta_2=\sum_i\beta_i\omega_i,\quad S\omega_i=\sum_ja_{ij}\omega_j.
\] 
Then the operators  $\omega\odot\eta_2$ and $\ast\iota(\omega)$, defined in the first section, are represented by matrices
\begin{align} 
\omega\odot\eta_2 &=
\left[\begin{array}{ccc} 2\alpha_1\beta_1 & \alpha_1\beta_2+\alpha_2\beta_1 & \alpha_1\beta_3+\alpha_3\beta_1 \\
\alpha_1\beta_2+\alpha_2\beta_1 & 2\alpha_2\beta_2 & \alpha_2\beta_3+\alpha_3\beta_2 \\
\alpha_1\beta_3+\alpha_3\beta_1 & \alpha_2\beta_3+\alpha_3\beta_2 & 2\alpha_3\beta_3 \end{array}\right], \label{mex1} \\
\ast\iota(\omega) &=\left[\begin{array}{ccc} 0 & -\alpha_3 & \alpha_2 \\ \alpha_3 & 0 & -\alpha_1 \\ -\alpha_2 & \alpha_1 & 0 \end{array}\right]. \label{mex2}
\end{align}
Put
\begin{align} \label{setsm}
U_1 &=\{x\in U: (a_{23}(x)-\mu(x)\beta_1(x))^2+(a_{13}(x)+\mu(x)\beta_2(x))^2>0 \}, \notag \\
U_2 &=\{x\in U: (a_{23}(x)+\mu(x)\beta_1(x))^2+(a_{12}(x)-\mu(x)\beta_3(x))^2>0 \}, \\
U_3 &=\{x\in U: (a_{12}(x)+\mu(x)\beta_3(x))^2+(a_{13}(x)-\mu(x)\beta_2(x))^2>0 \}. \notag
\end{align}
Then $U_1\cup U_2\cup U_3=U$. 

For two $1$--forms $\sigma$ and $\tau$ we write $\sigma\equiv \tau$ if there is nowhere vanishing smooth function $f$ such that $\sigma=f\tau$. We have
\begin{lem}\label{normalfield}
If a diffeomorphism $\varphi$ is conformal on a distribution $D$, then
\begin{align*}
\omega &\equiv(a_{13}+\mu\beta_2)\omega_1+(a_{23}-\mu\beta_1)\omega_2+(a_{33}-\lambda)\omega_3\quad\textrm{on $U_1$}, \\
\omega &\equiv(a_{12}-\mu\beta_3)\omega_1+(a_{22}-\lambda)\omega_2+(a_{23}+\mu\beta_1)\omega_3\quad\textrm{on $U_2$}, \\
\omega &\equiv(a_{11}-\lambda)\omega_1+(a_{12}+\mu\beta_3)\omega_2+(a_{13}-\mu\beta_2)\omega_3\quad\textrm{on $U_3$}.
\end{align*}
\end{lem}
\begin{proof} Proof is elementary but requires a lot of calculations. Details are left to the reader. By Theorem \ref{t1} $\varphi$ is conformal on $D$ if and only if
\begin{equation}\label{lc1}
B(\omega)=\mu(\omega\odot\eta_2).
\end{equation}
Moreover, 
\begin{equation}\label{lc2}
S\eta_2=\lambda\eta_2.
\end{equation}
Thus, using \eqref{mex1} and \eqref{mex2}, we get
\begin{align*}
(a_{23}-\mu\beta_1)\alpha_1 &=(a_{13}+\mu\beta_2)\alpha_2, \\
(a_{23}+\mu\beta_1)\alpha_1 &=(a_{12}-\mu\beta_3)\alpha_3, \\
(a_{12}+\mu\beta_3)\alpha_3 &=(a_{13}-\mu\beta_2)\alpha_2.
\end{align*}
Hence
\begin{align}
\alpha_1 &=C(a_{13}+\mu\beta_2), &\alpha_2 &=C(a_{23}-\mu\beta_1) & \textrm{on $U_1$}, \label{ln1} \\ 
\alpha_1 &=C'(a_{12}-\mu\beta_3), &\alpha_3 &=C'(a_{23}+\mu\beta_1) & \textrm{on $U_2$}, \label{ln2} \\
\alpha_2 &=C''(a_{12}+\mu\beta_3), &\alpha_3 &=C''(a_{13}-\mu\beta_2) & \textrm{on $U_3$}, \label{ln3}
\end{align}
for some $C, C'$ and $C''$. By \eqref{lc1} one can see that $C, C'$ and $C''$ are nowhere vanishing. Consider condition \eqref{ln1}. Since $\skal{\omega}{\eta_2}=0$,
\[
Ca_{13}\beta_1+Ca_{23}\beta_2+\alpha_3\beta_3=0
\] 
Moreover, by \eqref{lc2}
\[
Ca_{13}\beta_1+Ca_{23}\beta_2+C(a_{33}-\lambda)\beta_3=0.
\]
Hence
\[
\beta_3=0\quad\textrm{or}\quad \alpha_3=C(a_{33}-\lambda).
\]
Assuming $\beta_3=0$ and $\alpha_3 \not=C(a_{33}-\lambda)$ and using \eqref{lc1} and \eqref{lc2}, after some calculations we get a contradiction. Finally
\[
\alpha_1=C(a_{13}+\mu\beta_2),\quad \alpha_2=C(a_{23}-\mu\beta_1),\quad \alpha_3=C(a_{33}-\lambda) \quad \textrm{on $U_1$}.
\]
Analogously we consider conditions \eqref{ln2} and \eqref{ln3}.
\end{proof}
Let
\[
d a_{ij}=\sum_ka^k_{ij}\omega_k, \quad d\mu=\sum_k\mu_k\omega_k, \quad d\beta_i=\sum_k\beta_{ik}\omega_k, \quad d\lambda=\sum_k\gamma_k\omega_k
\]
Let $[e_j,e_k]=\sum_iC^i_{jk}e_i$. Then
\[
d\omega_i=-\sum_{j<k}C^i_{jk}\omega_j\wedge\omega_k. 
\]
Therefore by Lemma \ref{normalfield} integrability condition $d\omega\wedge\omega=0$ is on $U_1$, $U_2$ and $U_3$ respectively
\begin{gather} 
\begin{split} \label{integr1}
0 &=\bigg(-a^2_{13}+a^1_{23}-\beta_1\mu_1-\beta_2\mu_2-\mu\beta_{11}-\mu\beta_{22} \\
  &-(a_{13}+\mu\beta_2)C^1_{12}-(a_{23}-\mu\beta_1)C^2_{12}-(a_{33}-\lambda)C^3_{12} \bigg)(a_{33}-\lambda) \\
  &-\bigg(-a^3_{13}+a^1_{33}-\beta_2\mu_3-\mu\beta_{23}-\gamma_1 \\
  &-(a_{13}+\mu\beta_2)C^1_{13}-(a_{23}-\mu\beta_1)C^2_{13}-(a_{33}-\lambda)C^3_{13} \bigg)(a_{23}-\mu\beta_1) \\
  &+\bigg(-a^3_{23}+a^2_{33}+\beta_1\mu_3+\mu\beta_{13}-\gamma_2 \\
  &-(a_{13}+\mu\beta_2)C^1_{23}-(a_{23}-\mu\beta_1)C^2_{23}-(a_{33}-\lambda)C^3_{23} \bigg)(a_{13}+\mu\beta_2), \\
\end{split} \\
\begin{split} \label{integr2}
0 &=\bigg(-a^2_{12}+a^1_{22}+\beta_3\mu_2+\mu\beta_{32}-\gamma_1 \\
  &-(a_{12}-\mu\beta_3)C^1_{12}-(a_{22}-\lambda)C^2_{12}-(a_{23}+\mu\beta_1)C^3_{12} \bigg)(a_{23}+\mu\beta_1) \\
  &-\bigg(-a^3_{12}+a^1_{23}+\beta_3\mu_3+\beta_1\mu_1+\mu\beta_{33}+\mu\beta_{11} \\
  &-(a_{12}-\mu\beta_3)C^1_{13}-(a_{22}-\lambda)C^2_{13}-(a_{23}+\mu\beta_1)C^3_{13} \bigg)(a_{22}-\lambda) \\
  &+\bigg(-a^3_{22}+a^2_{23}+\beta_1\mu_2+\mu\beta_{12}+\gamma_3 \\
  &-(a_{12}-\mu\beta_3)C^1_{23}-(a_{22}-\lambda)C^2_{23}-(a_{23}+\mu\beta_1)C^3_{23} \bigg)(a_{12}-\mu\beta_1), \\
\end{split} \\
\begin{split} \label{integr3}
0 &=\bigg(-a^2_{11}+a^1_{12}+\beta_3\mu_1+\mu\beta_{31}+\gamma_2 \\
  &-(a_{11}-\lambda)C^1_{12}-(a_{12}+\mu\beta_3)C^2_{12}-(a_{13}-\mu\beta_2)C^3_{12} \bigg)(a_{13}-\mu\beta_2) \\
  &-\bigg(-a^3_{11}+a^1_{13}-\beta_2\mu_1-\mu\beta_{21}+\gamma_3 \\
  &-(a_{11}-\lambda)C^1_{13}-(a_{12}+\mu\beta_3)C^2_{13}-(a_{13}-\mu\beta_2)C^3_{13} \bigg)(a_{12}+\mu\beta_3) \\
  &+\bigg(-a^3_{12}+a^2_{13}-\beta_2\mu_2-\beta_3\mu_3-\mu\beta_{22}-\mu\beta_{33} \\
  &-(a_{11}-\lambda)C^1_{12}-(a_{23}+\mu\beta_3)C^2_{23}-(a_{13}-\mu\beta_2)C^3_{23} \bigg)(a_{11}-\lambda).
\end{split} 
\end{gather}
Thus above considerations and Theorem \ref{t1} imply
\begin{thm}\label{t2} Let $\varphi:M\to N$ be a diffeomorphism between $3$--dimensional oriented Riemannian manifolds. Suppose $S=(\varphi_{\ast})^*\varphi_{\ast}$ has distinct eigenvalues $\lambda_1<\lambda_2<\lambda_3$. Let $\xi_2$ be the unit eigenvector corresponding to $\lambda_2$ and let $\eta_2$ be a $1$--form dual to $\xi_2$. If $\varphi$ is leafwise conformal on an open subset $U$ of $M$ then conditions \eqref{integr1}--\eqref{integr3} hold, where $U_1, U_2, U_3$ are defined by \eqref{setsm} and $\mu$ is given by \eqref{mu} with the sign $+$ or $-$ instead of $\pm$. 
\end{thm}

\section{Some necessary and sufficient conditions of integrability}
Let $\varphi:M\to N$ be a diffeomorphism between Riemannian manifolds. Let $\lambda_1,\lambda_2,\lambda_3$ be the eigenvalues of the operator $S=(\varphi_{\ast})^*\varphi_{\ast}:TM\to TM$ and $\xi_1,\xi_2,\xi_3$ be the corresponding unit eigenvectors. Assume $\lambda_1<\lambda_2<\lambda_3$. Let $\eta_1,\eta_2,\eta_3$ be the basis dual to $\xi_1,\xi_2,\xi_3$. Assume $\xi_i$ and $\eta_i$ are globally smooth. Consider $1$--forms $\omega_{\pm}$ given by \eqref{eq1} and put $D_{\pm}={\rm ker}\omega_{\pm}$. Then by Lemma \ref{l1} $\varphi$ is conformal on the distributions $D_{\pm}$. We study the integrability condition $\omega_{\pm}\wedge d\omega_{\pm}=0$. After simple calculations we get
\begin{equation}\label{pr1}
\eta_1\wedge d\eta_1+\chi^2\eta_3\wedge d\eta_3=\pm d(\chi\eta_1\wedge\eta_3 ),
\end{equation}
where
\[
\chi=\frac{\sqrt{\lambda_3-\lambda_2}}{\sqrt{\lambda_2-\lambda_1}}.
\]
Therefore we have
\begin{prop}\label{p2} 
If $D_{\pm}$ are both integrable, then $\eta_1\wedge d\eta_1+\chi^2\eta_3\wedge d\eta_3=0$ and $2$--form $\chi\eta_1\wedge\eta_3$ is closed.
\end{prop}
Write $\eta_1$ and $\eta_3$ in terms of $\omega_{\pm}$,
\[
\eta_1=\frac{1}{2}\frac{\sqrt{\lambda_3-\lambda_1}}{\sqrt{\lambda_2-\lambda_1}}(\omega_++\omega_-),\qquad
\eta_3=\frac{1}{2}\frac{\sqrt{\lambda_3-\lambda_1}}{\sqrt{\lambda_3-\lambda_2}}(\omega_+-\omega_-).
\]
Then
\begin{equation} \label{clform}
\chi\eta_1\wedge\eta_3=-\frac{1}{2}(1+\chi^2)\omega_+\wedge\omega_-,
\end{equation}
\begin{prop} \label{p3}
Assume $M$ is orientable and closed. Suppose $\eta_i$, $\chi$, $D_{\pm}$ are smooth and globally defined on $M$. If both $D_{\pm}$ are integrable then
\[
\int_M\xi_2({\rm log}\chi) {\rm vol}=0,
\]
where ${\rm vol}$ is the volume element on $M$.
\end{prop}
\begin{proof} By Proposition \ref{p2}, $d(\chi\eta_1\wedge\eta_3)=0$. Hence
\[
\frac{1}{\chi}(\xi_2\chi) \eta_1\wedge\eta_2\wedge\eta_3=d(\eta_1\wedge\eta_3)
\]
\end{proof}

\begin{ex} \label{e2}Consider a diffeomorphism $\varphi:\mathbb{R}^3\to\mathbb{R}^3$ such that $S$ is diagonal in the canonical basis $e_1,e_2,e_3$, $S={\rm diag}(\lambda_1,\lambda_2,\lambda_3)$. Assume $\lambda_1<\lambda_2<\lambda_3$. Then $\eta_i=dx_i$, $i=1,2,3$. Therefore, the integrability condition reduces to
\[
\frac{\partial\chi}{\partial x_2}=0,
\]
which we can write in the form
\[
(\lambda_2-\lambda_1)(\frac{\partial\lambda_3}{\partial x_2}-\frac{\partial\lambda_2}{\partial x_2})=
(\lambda_3-\lambda_2)(\frac{\partial\lambda_2}{\partial x_2}-\frac{\partial\lambda_1}{\partial x_2}).
\]
Since $\lambda_i=\sum_j(\partial\phi_j/\partial x_i)^2$ we obtain
\begin{multline*}
\sum_j((\frac{\partial\phi_j}{\partial x_2})^2-(\frac{\partial\phi_j}{\partial x_1})^2) 
\sum_k(\frac{\partial\phi_j}{\partial x_3}\frac{\partial^2\phi_j}{\partial x_3\partial x_2}-\frac{\partial\phi_j}{\partial x_2}\frac{\partial^2\phi_j}{\partial x_2^2})
= \\
\sum_k(\frac{\partial\phi_j}{\partial x_2}\frac{\partial^2\phi_j}{\partial x_2^2}-\frac{\partial\phi_j}{\partial x_3}\frac{\partial^2\phi_j}{\partial x_1\partial x_2})
\sum_j((\frac{\partial\phi_j}{\partial x_3})^2-(\frac{\partial\phi_j}{\partial x_1})^2).
\end{multline*}
\end{ex}

\begin{ex} \label{e3} Consider a local diffeomorphism from Example \ref{e1}. Then $\chi=1$ and by (\ref{clform}) we have
\[
\chi \eta_1\wedge\eta_3=\frac{\sqrt{2}}{\sqrt{2+2\sin^2(x_2+x_3)}}dx_1\wedge dx_2-\frac{\sqrt{2}\sin(x_2+x_3)}{\sqrt{2+2\sin^2(x_2+x_3)}}dx_2\wedge dx_3.
\]
Thus the above form is not closed. By Proposition \ref{p2} one of the distributions $D_{\pm}$ is not integrable. Since $D_+$ is obviously integrable, we get that $D_-$ is not integrable.
\end{ex}

\section{Local leafwise holomorphicity}
A coordinate system $\psi$ on an $n$--dimensional Riemannian manifold $(M,g)$ is called {\it foliated conformal chart} if the map $z\mapsto\psi^{-1}(z,q)$, $z\in\mathbb{R}^2$, $q\in\mathbb{R}^{n-2}$, is conformal for all $q$.

The aim of this section is to prove the following result.
\begin{thm} \label{t3} A map $\varphi:M\to N$ between Riemannian manifolds is locally leafwise conformal if for every $x\in M$ there are foliated conformal charts $\psi$ and $\tilde{\psi}$ in neighbourhoods of $x$ and $\varphi(x)$ respectively, such that
\[
\tilde{\psi}\circ \varphi\circ \psi^{-1}(z,q)=(h(z,q),q),
\]
where for every $q\in\mathbb{R}^{n-2}$ the map $\mathbb{R}^2\ni z\mapsto h(z,q)\in\mathbb{R}^2$ is holomorphic.
\end{thm}

Let us first review some facts about the Beltrami equation and isothermal coordinates. 

Assume all considered functions are smooth. By the {\it Beltrami equation} we mean the equation
\[
\frac{\partial w}{\partial \bar z}=\mu\frac{\partial w}{\partial z},
\]
where $\mu,w:\mathbb{C}\to\mathbb{C}$. If $|\mu|<k<1$ for some $k$, then the Beltrami equation has a unique smooth solution $w^{\mu}$ which leaves $0$, $1$ and $\infty$ fixed. Moreover $w^{\mu}$ has positive Jacobian, see \cite{bjs}. We have also smooth dependence of solutions of the Beltrami equation \cite{es}.
\begin{thm}[Riemann's mapping theorem for variable metric] \label{t4} For each positive $k<1$ the map $\mu\mapsto w^{\mu}$ is a homeomorphism of the set $\{\mu\in C^{\infty}(\mathbb{C},\mathbb{C}): \sup|\mu|<k\}$ onto its image in $C^{\infty}(\mathbb{C},\mathbb{C})$.
In particular, the map
\[
w:\mathbb{C}\times T\ni (z,t)\mapsto (w^{\mu_t}(z),t)\in \mathbb{C}\times T
\]
is a diffeomorphism, where $T$ is open subset of $\mathbb{R}^q$, $\mu_t(z)=\mu(z,t)$.
\end{thm}
In Theorem \ref{t4}, $C^{\infty}(\mathbb{C},\mathbb{C})$ is a Frechet space of all smooth functions $f:\mathbb{C}\to\mathbb{C}$ with $C^{\infty}$ topology.

Consider $\mathbb{R}^2=\mathbb{C}$ with a Riemannian metric $g=Edx^2+2Fdxdy+Gdy^2$, where $E>0$, $EG-F^2>0$. A coordinate system $w=(u,v)$ is called \emph{isothermal} if there is a positive function $\lambda$ such that
\[
g=\lambda(du^2+dv^2).
\] 
Put
\begin{equation} \label{rem1}
\mu=\frac{E-G+2iF}{E+G+2\sqrt{EG-F^2}}.
\end{equation}
Take a closed ball $K\subset\mathbb{C}$. Since $|\mu|<1$, $\sup_K|\mu|<k_0<1$ for some $k_0$. Extend $\mu$ smoothly to the whole plane in such a way that $\sup|\mu|<k<1$ for some $k$. Then $w=w^{\mu}|\textrm{int}K$ is an isotermal coordinate system for $g$, see \cite{sp}. For a foliation by planes $L_t=\mathbb{C}\times\{t\}, t\in T$, with a Riemannian metric $g$, on each leaf $L_t$ we have $g|_{L_t\times L_t}=E_tdx^2+2F_tdxdy+G_tdy^2$. Therefore, in the same way as before, by Theorem \ref{t4} for $\mu_t$ defined by (\ref{rem1}), there is a coordinate system $w_t(z)=w(z,t)$ such that $w_t:L_t\to L_t$ is isothermal for every $t$. We say that $w$ is a \emph{foliated isothermal coordinate system}.

\begin{proof}[Proof of Theorem \ref{t3}]
Let $\varphi:M\to N$ be an $\mathcal{F}$--conformal diffeomorphism, $(M,g_M)$ and $(N,g_N)$ Riemannian manifolds, $\mathcal{F}$ a  $2$--dimensional foliation on $M$. Fix $x\in M$ and put $y=\varphi(x)$. Let $\chi$ be a foliated map in a neighbourhood of $x$ and let $\rho=\chi\circ \varphi^{-1}$. Then $\rho\circ \varphi\circ\chi^{-1}$ is the identity map on the foliation
\[
\mathcal{F}_0=\chi(\mathcal{F})=\{U\times\{t\}\}_{t\in T},
\]
where $U$ is an open subset of $\mathbb{C}$ and $T$ open subset of $\mathbb{R}^{{\rm codim}\mathcal{F}}$. Obviously, ${\rm id}=\rho\circ \varphi\circ\chi^{-1}$ is a $\mathcal{F}_0$--conformal map with respect to Riemannian metrics $(\chi^{-1})^*g_M$ and $(\rho^{-1})^*g_N$. Let $w_M$ and $w_N$ be foliated isothermal coordinate systems on $U$, shrinking $U$ if necessary, for $(\chi^{-1})^*g_M$ and $(\rho^{-1})^*g_N$ respectively. Then
\[
h=w_N\circ\rho\circ \varphi\circ\chi^{-1}\circ w_M^{-1}=w_N\circ w_M^{-1} 
\]
is a $\mathcal{F}_0$--conformal diffeomorphism with respect to Riemannain metrics $(\chi^{-1}\circ w_M^{-1})^*g_M$ and $(\rho^{-1}\circ w_N^{-1})^*g_N$, which on leaves of $\mathcal{F}_0$ are conformal with Euclidean metric. Thus $h$ is $\mathcal{F}_0$--conformal with respect to Euclidean metric. Therefore, maps $h:L_t\to L_t$, $L_t=U\times\{t\}$, $t\in T$, are all holomorphic or all antiholomorphic. If $h:L_t\to L_t$, $t\in T$, are antiholomorphic, we replace $h$ by
\[
\tilde{h}=\tau\circ h,
\]
where $\tau(z,t)=(\bar{z},t)$. Then $\psi=w_M\circ\chi$ and $\tilde{\psi}=w_N\circ\rho$ are desired.
\end{proof}

\begin{ackn} This article is based on a part of author's PhD Thesis \cite{kn2}. The author wishes to thank his advisor Professor Antoni Pierzchalski for helpful discussions.
\end{ackn}

\end{document}